\documentclass[11pt]{article}
\usepackage{}
\usepackage{amsfonts}
\usepackage{amssymb}
\usepackage{amsmath,amssymb,graphicx}
\usepackage{caption2}

\usepackage{epsfig}
\usepackage{pstricks}
\usepackage{pst-node}
\usepackage{pst-tree}
\usepackage{pst-plot}
\usepackage{pst-text}
\usepackage{multido}
\usepackage{subfigure}
\newtheorem{theorem}{Theorem}[section]
\newtheorem{lemma}[theorem]{Lemma}

\newtheorem{corollary}[theorem]{Corollary}

\newcommand{\be}{\begin{eqnarray}}
\newcommand{\ee}{\end{eqnarray}}
\newcommand{\ba}{\begin{array}}
\newcommand{\ea}{\end{array}}
\newcommand{\ben}{\begin{eqnarray*}}
\newcommand{\een}{\end{eqnarray*}}
\newcommand{\mc}{\mathcal}
\newcommand{\mb}{\mathbb}
\newcommand{\nn}{\nonumber}
\newcommand{\sn}{\mathrm{sn}}
\newcommand{\cn}{\mathrm{cn}}
\newcommand{\dn}{\mathrm{dn}}
\newcommand{\tr}{\mathrm{tr}}
\newcommand{\ls}{\;\!}
\newcommand{\llog}{\mathrm{Log}}
\newcommand{\K}{\mathcal{K}}
\newcommand{\sq}{\sqrt}

\definecolor{hylightgray}{rgb}{0.9,0.9,0.9}

\begin{document}

\title{{\large \bf Uniformisation of a once-punctured annulus}
}
\author{\normalsize  Tanran Zhang }
\date{Graduate School of Information Sciences, Tohoku University\\
      Aoba-ku, Sendai 980-8579, Japan\\
      \textsf{zhang@ims.is.tohoku.ac.jp}\\}
\maketitle \baselineskip 21pt
\noindent

\begin{minipage}{138mm}
\renewcommand{\baselinestretch}{1} \normalsize
\textbf{Abstract.} The universal cover or the covering group of a hyperbolic Riemann surface $X$ is important but hard to express explicitly. It can be, however, detected by the uniformisation and a suitable description of $X$. Beardon proposed five different ways to describe twice-punctured disks using fundamental domain, hyperbolic length, collar and extremal length in 2012. We parameterize a once-punctured annulus $A$ in terms of five parameter pairs and give explicit formulas about the hyperbolic structure and the complex structure of $A$. Several degenerating cases are also treated.
\vspace*{2mm}

\vspace*{2mm}

\end{minipage}\\
\renewcommand{\baselinestretch}{1} \normalsize
\section{Introduction}
The Uniformisation Theorem implies that every Riemann surface $X$ is conformal equivalent to the quotient space $\mb{H}/G$, where $G$ is a torsion-free Fuchsian group acting on the upper half plane $\mb{H}:=\{z \in \mb{C}: \textrm{Im} z >0\}$, if $X$ is not conformally equivalent to the Riemann sphere $\hat{\mb{C}}$, the complex plane $\mb{C}$, the once-punctured complex plane $\mb{C}\backslash \{a\}$ or a complex torus. It is, however, difficult to find an explicit form of the holomorphic universal cover $\pi$ or the covering group $G$, except for several special cases (see e.g. \cite{Hejhal1974, Sugawa1998uni}). For a twice-punctured unit disk, Hempel and Smith \cite{Hempelsmith1988, Hempel1989acc, Hempel1989uni} considered the uniformisation problem and the hyperbolic metric, and Beardon \cite{Beardon2012} provided five parameters to characterize the twice-punctured disk via the hyperbolic structure and the complex structure of it. Nevanlinna \cite[I.3, I.4]{Nevanlinna1970} introduced a method to regard the puncture as the extremal case when a boundary curve shrinks to a single point. In this article we give five parameter pairs to uniformize a once-punctured annulus $A$. These parameter pairs can be divided into two classes which are corresponding to the hyperbolic structure and complex structure of $A$, respectively.

Let $\gamma$ be a simple closed geodesic on a hyperbolic surface $X$, and let
\be \label{collar}C_{\theta}(\gamma):=\{x\in X: \delta_{X}(x,\gamma) < \sinh^{-1}(\tan \theta) /2\},
\ee
where $\delta_{X}$ is the hyperbolic distance on $X$ of the Gaussian curvature $-1$. $C_{\theta}(\gamma)$ is called a $collar$ \emph{about} $\gamma$ \emph{of angular width} $\theta$ if it is doubly connected. When $X \backslash \gamma$ has a doubly connected component $W$ and $\gamma$ is homotopic to no puncture, $\gamma$ is homotopic to a border of $X$, and $\gamma$ is called $peripheral$. If $C_{\theta}(\gamma)$ is a collar, $\widetilde{C}_{\theta}(\gamma):=C_{\theta}(\gamma) \cup W$ is a doubly connected subdomain in $X$ containing $\gamma$. We will refer to $\widetilde{C}_{\theta}(\gamma)$ as a \emph{peripheral collar about $\gamma$ of angular width $\theta$}. In Section 2 we will give the details of peripheral collars.

We will use Legendre's complete elliptic integrals $K(r)$ of the first kind. Let $K'(r):=K(r')=K(\sqrt{1-r^2})$ for $0<r<1$, and $\mu(r):={\pi K'(r) } / {\left(2 K(r) \right)}$. The properties of $K(r)$, $K'(r)$ and $\mu(r)$ are given in Section 2.

After an application of a rotation and a similarity, we only need to consider the punctured annulus
\be \label{natural A}
A:=\{z:1/R <|z|<R\}\backslash \{a\},\ R>1,\ 1/R<a<R.
\ee
We denote $C_1:=\{z:|z|=1/R\}$,\ $C_2:=\{z:|z|=R\}$, and let $\mc{C}_1$, $\mc{C}_2$ be the free homotopy classes of the circles $\{z: |z|=r_1\}$, $\{z: |z|=r_2\}$ in $A$, respectively, where $a<r_1<R$, $1/R<r_2<a$. So $\mc{C}_1$ separates $C_1 \cup \{a\}$ from $C_2$, and $\mc{C}_2$ separates $C_2 \cup \{a\}$ from $C_1$. Let $\gamma_1$, $\gamma_2$ be the hyperbolic geodesics in $\mc{C}_1$, $\mc{C}_2$. The main results in this article are as follows.

\begin{theorem}\label{main theorem}
\textsl{Let $l_1$, $l_2$ be the hyperbolic lengths of $\gamma_1$, $\gamma_2$, and $\theta_1$, $\theta_2$ be the angular widths of the maximal peripheral collars about $\gamma_1$, $\gamma_2$ in the punctured annulus $A$. Then $(\l_1, l_2)$ and $(\theta_1, \theta_2)$ satisfy
\be \label{main quantity}
\displaystyle \cos \theta_1=\frac{\sinh \frac{l_1}{2}} {\cosh \frac{l_1}{2}+\cosh \frac{l_2}{2}},\ \ \ \cos \theta_2=\frac{\sinh \frac{l_2}{2}} {\cosh \frac{l_1}{2}+\cosh \frac{l_2}{2}}.
\ee}
\end{theorem}

\begin{theorem} \label{extremal length theorem}
\textsl{Let $\lambda_1$ and $\lambda_2$ be the extremal lengths of $\mc{C}_1$ and $\mc{C}_2$. Select a positive number $q$ such that $\mu(q)=4\log R$ and let
$\K:=K(q)$, $\K ':=K'(q)$.
Then
\be \label{extremal length results}
\lambda_1=\frac{2\pi}{\mu(p_1)},\ \ \ \lambda_2=\frac{2\pi}{\mu(p_2)},
\ee
where
\be \label{p1p2}
p_1=\frac{\sq q(\dn\,u_1+1)}{q+\dn\,u_1}, \quad p_2=\frac{\sq q(\dn\,u_2+1)}{q+\dn\,u_2}
\ee
with
\be \label{u1 u2}
u_1=\frac{2 \K}{\pi}\log Ra,\ \ \  u_2=\frac{2 \K}{\pi}\log \frac{R}{a},
\ee
and the Jacobian elliptic function $\dn$ in \eqref{p1p2} has the modulus ${q}'=\sqrt{1-q^2}$.}
\end{theorem}

Theorems \ref{main theorem} and \ref{extremal length theorem} reveal the connections between $(\theta_1, \theta_2)$ and $(l_1,l_2)$, $(\lambda_1, \lambda_2)$ and $(R,a)$, respectively. This implies that we have two classes of parameter pairs to describe $A$, one class related to $(l_1,l_2)$ and the other to $(R,a)$. Actually the two parameter classes are corresponding to two kinds of structures on $A$, hyperbolic and complex structures.

This article is organized in the following way. Section 2 is about preliminaries, and we construct the covering group $G$ and the fundamental domain $D_A$ of $A$ on $\mathbb{H}$ there. The construction of $G$ gives rise to two real numbers $r$ and $k$, which also form a parameter pair of $A$. Section 3 is devoted to the hyperbolic structure of $A$. We prove Theorem \ref{main theorem} using the parameter pair of the covering group $G$ given in Section 3. Section 4 is devoted to the complex structure of $A$, and we prove Theorem \ref{extremal length theorem} there. Section 5 is about the degenerating cases when the puncture is tending to one of the boundary circles, and when the boundaries are shrinking to points. We give an observation for the once-punctured annulus $A$ when only one boundary circle of $A$ is shrinking to a single point, in which case $A$ is becoming a twice-punctured disk.

\section{Preliminary}
Let $\mb{D}:=\{z \in \mb{C}: |z|<1\}$. The hyperbolic metrics on $\mb{H}$ and $\mb{D}$ of the Gaussian curvature $-1$ are given by
\be
\rho_{\mb{H}}(z)|dz|=\frac {|dz|}{\textrm{Im} z},\quad \rho_{\mb{D}}(z)|dz|=\frac {2|dz|}{1-{|z|}^2}.
\ee
They induce the hyperbolic lengths $\delta_{\mb{H}}(x,y)$ and $\delta_{\mb{D}}(x,y)$ between two points (or sets) $x$ and $y$ in $\mb{H}$ and $\mb{D}$. The universal covering space $\widetilde{X}$ of a Riemann surface $X$ can be tesselated by a fundamental domain and its images under the covering group $G$ acting on $\widetilde{X}$. For a hyperbolic surface $X$, we take $\widetilde{X}$ to be the upper half plane $\mathbb{H}$ or the open unit disk $\mathbb{D}$. A domain $D \subseteq \widetilde{X}$ is called a \emph{fundamental domain} for $G$ if $D$ satisfies the following two conditions: $g(D) \cap D =\emptyset$ for all $g \in G$, $g\neq \textrm{Id}$; $\bigcup_{g \in G}\overline{g(D)}=\widetilde{X}$. The hyperbolic metric on $\widetilde{X}$ can be projected under the quotient mapping to a metric on $X$ which is called the \emph{hyperbolic metric on $X$}. It is the intrinsic metric on $X$. It is independent of the choice of the universal cover from $\widetilde{X}$.

We identify the M\"{o}bius transformation
\be \label{mobius trans}
\phi(z)=\frac{az+b}{cz+d},\ \ ad-bc=1, \ \ a,b,c,d \in \mb{C},
\ee
with the $2\times 2$ complex matrices
$\pm \left(
\begin{array}{cc}
a& b  \\
c& d
\end{array}
\right)\in{\textrm{PSL}(2,\mb{C})}$
which is also denoted by $\phi$, and define the \emph{trace} of $\phi$ by $\tr\,\phi=\pm (a+d)$, such that $\tr^2 \phi=(a+d)^2$ is well defined.  M\"{o}bius transformations preserve the hyperbolic metric. All the conformal isometries of $\mb{H}$ are M\"{o}bius transformations with $a,\;b,\;c,\;d$ being real numbers. The \emph{translation length} of $\phi$ is defined by $T(\phi)=\inf_{z\in \mb{H}} \delta_{\mb{H}}(z,\phi(z))$, where $\delta_{\mb{H}}$ is the hyperbolic distance on $\mb{H}$. When $\phi$ is hyperbolic, $T(\phi)>0$. It is known that $2\cosh(T(\phi)/2)=|\tr \,\phi|$ (see \cite[7.34]{Beardonbookdiscrete}).

For a simple closed geodesic $\gamma$ of hyperbolic length $l$ in a hyperbolic surface $X$, let $C_{\theta}(\gamma)$ be a collar about $\gamma$ of width $\theta$, $0<\theta < \frac{\pi}{2}$. Then there exists a hyperbolic transformation $f$ in the covering group $G$ of $X$, such that $f(C_{\theta}(\gamma))=C_{\theta}(\gamma)$ and $\mathrm{Stab}_G \left(C_{\theta}(\gamma)\right)=\langle f \rangle$, where $\mathrm{Stab}_G \left(C_{\theta}(\gamma)\right)$ is the stabilizer of $C_{\theta}(\gamma)$ in $G$.  Moreover, $C_{\theta}(\gamma)$ satisfies $h(C_{\theta}(\gamma)) \cap C_{\theta}(\gamma)=\emptyset$ if $h \in G \backslash \langle f \rangle$. By conjugation we may assume that $f: z \mapsto k^2 z$, $k=\exp {\frac{l}{2}}$, then $\{z: 1<|z|<k^2, \frac{\pi}{2}-\theta < \textrm{arg}\; z <\frac{\pi}{2}+\theta\}$ is a fundamental domain for $\langle f \rangle$ in $C_{\theta}(\gamma)$. The collar lemma showed that if $\theta$ satisfies
\be \label{kl collar}
\tan \theta \leq \frac{2}{k-k^{-1}}
\ee
$C_{\theta}$ is a collar about $\gamma$ of width $\theta$ (see \cite[Lemma 7.7.1]{Keenbook}). In addition, if $\gamma$ is peripheral, $\{z: 1<|z|<k^2, \frac{\pi}{2}-\theta < \textrm{arg}\; z < \pi \}$ (or $\{z: 1<|z|<k^2, 0 < \textrm{arg}\; z < \pi -\theta\}$ by conjugation) is a fundamental domain for $\langle f \rangle$ in the peripheral collar $\widetilde{C}_{\theta}(\gamma)$ which contains $C_{\theta}(\gamma)$.

We have the following lemma for the covering group $G$ of $A$.
\begin{lemma}\label{fundamental domain lemma}
\textsl{For the punctured annulus $A$ given by \eqref{natural A}, there exist two real numbers $k$ and $r$, $1<r<k$, such that the group $G$ generated by the hyperbolic transformation $f$ and the parabolic transformation $g$ is the covering group of $A$ acting on $\mb{H}$, where
\be \label{transformations}
f=\pm \left(
\begin{array}{cc}
k& 0  \\
0& k^{-1}
\end{array}
\right), \ \ \ g=\frac{\pm 1}{r-1}\left(
\begin{array}{cc}
2r& -(r+1)  \\
r+1& -2
\end{array}
\right).
\ee
}
\end{lemma}
\textbf{Proof.} We first construct a fundamental domain of the covering group acting on $\mb{D}$ and a universal cover $\pi: \mb{D}\rightarrow A$ of the punctured annulus $A$. $A$ can be divided into two pieces, the upper half $A^+$ and the lower half $A^-$, by the three Euclidean line segments $\eta_1:=(-R, -1/R)$, $\eta_2:=(1/R, a)$, $\eta_3:=(a, R)$ lying on the real axis. Since $A$ is symmetric with respect to the real axis, each lift of $\eta_i$ under $\pi$ is a hyperbolic line in $\mb{D}$ for $i=1,2,3$ (see Hemple and Smith \cite[Section 2]{Hempel1989uni}). We take one component $D^0_{A^+}$ of the pre-image of $A^+$ under $\pi$, then the inverse of $\pi$ has a univalent branch $\tau: A^+\rightarrow D^0_{A^+}$. We now describe $\partial D^0_{A^+} \cap \mb{D}$. The function $\tau$ can be extended continuously to $\overline{A^+}\cap A=A^+\cup \eta_1 \cup \eta_2 \cup \eta_3$. We denote the extension of $\tau$ still by $\tau$, so that $\tau(\eta_1)\cup \tau(\eta_2)\cup \tau(\eta_3)=\partial D^0_{A^+} \cap \mb{D}$. Since $\tau(\eta_i)$ is a hyperbolic line in $\mb{D}$ for $i=1,2,3$, by conjugacy of M\"{o}bius transformations if necessary, we may assume that $\tau(\eta_3)$ is the diameter $(-1,1)$ of $\mb{D}$, with $-1=\tau(a)$ and $1=\tau(R)$. Then $\tau(\eta_2)$ and $\tau(\eta_3)$ are hyperbolic lines in $\mb{D}$, where $\tau(\eta_2)$ and $\tau(\eta_3)$ are parallel with the common endpoint $-1$, and $\tau(\eta_1)$ is disjoint with $\tau(\eta_2)$ or $\tau(\eta_1)$, in accordance with the position of $\eta_1$, $\eta_2$ and $\eta_3$ in $A^+$. In such a way the upper half punctured annulus $A^+$ is mapped onto the region $D^{0}_{A^+}$ bounded by $\tau(\eta_i)$, $i=1,2,3$, and the two arcs on $\partial \mb{D}$ joining $1$ and $\tau(-R)$, $\tau(-1/R)$ and $\tau(1/R)$, shown in Figure \ref{fund domain in disk}.

\begin{figure}[htbp]
 \centering
 \subfigure[\label{fund domain in disk}]{\scalebox{0.7}[0.7]{
\begin{pspicture}(-3.2,-3.2)(3.2,3.5)
\pscircle[fillstyle=solid,fillcolor=hylightgray,linecolor=white](0,0){3.2}
\pscircle[fillstyle=solid,fillcolor=white,linecolor=white](-3.2,1.12){1.12}
\psarc[fillstyle=solid,fillcolor=white,linecolor=white](1.6,3.36){1.90}{170}{320}
\pscircle[fillstyle=solid,fillcolor=white,linecolor=white](-3.2,-1.12){1.12}
\pscircle[fillstyle=solid,fillcolor=white,linecolor=white](1.6,-3.36){1.90}

\pscircle(0,0){3.2}
\psline[linestyle=dashed](-3.2,0)(3.2,0)
\psarc(-3.2,1.12){1.12}{-90}{50}
\psarc(1.6,3.36){1.9}{187}{303}
\psarc(-3.2,-1.12){1.12}{-50}{90}
\psarc(1.6,-3.36){1.9}{57}{173}
\psarc[linestyle=dashed](6.56,0){5.6}{151}{209}
\psline(1.52,2.77)(1.46,2.64)(1.6,2.56)
\psline(1.52,-2.77)(1.46,-2.64)(1.6,-2.56)
\psline(1.24,1.7)(1.38,1.64)(1.33,1.5)
\psline(1.24,-1.7)(1.38,-1.64)(1.33,-1.5)

\uput[d](-3.55,0.2){$-1$}
\uput[d](4.06,0.26){$1=\tau(R)$}
\uput[d](3.26,2.32){$\tau(-R)$}
\uput[d](-0.19,3.9){$\tau(-\frac 1 R)$}
\uput[d](-3,2.5){$\tau(\frac 1 R)$}
\uput[d](-0.8,2.6){$D^0_{A^+}$}
\uput[d](-0.8,-1.12){$D^0_{A^-}$}
\uput[d](-1.6,1.3){$\tau(\eta_2)$}
\uput[d](0.52,1.7){$\tau(\eta_1)$}
\uput[d](0.48,-0.9){$\overline{\tau(\eta_1)}$}
\uput[d](-0.3,0.7){$\tau(\eta_3)$}
\uput[d](1.8,3.4){$\zeta_0$}
\uput[d](1.8,-2.7){$\bar{\zeta_0}$}
\end{pspicture}}  }\ \ \hspace{2cm}
\subfigure[\label{fund domain in upper}]{\scalebox{0.75}[0.75]{
\begin{pspicture}(-3.9,-0.64)(4.4,4)

\psarc[fillstyle=solid,fillcolor=hylightgray,linecolor=white](0,0){4}{0}{180}
\psarc[fillstyle=solid,fillcolor=white,linecolor=white](0,0){0.4}{0}{180}
\psarc[fillstyle=solid,fillcolor=white,linecolor=white](1.32,0){0.36}{0}{180}
\psarc[fillstyle=solid,fillcolor=white,linecolor=white](2.34,0){0.66}{0}{180}

\psline(-4.4,0)(4.4,0)
\psline{->}(0.33,0.44)(2.08,3.12)
\psarc(0,0){4}{0}{180}
\psarc[linestyle=dashed](0,0){1.68}{0}{180}
\psarcn{->}(2,0){0.8}{150}{80}
\psarc(0,0){0.4}{0}{180}
\psarc(1.32,0){0.36}{0}{180}
\psarc(2.34,0){0.66}{0}{180}

\uput[d](-4,-0.02){$-k$}
\uput[d](-0.42,0.02){$-\frac 1 k$}
\uput[d](0.39,0.02){$\frac 1 k$}
\uput[d](1.03,0.02){$\frac 1 r$}
\uput{6pt}[d](1.68,-0.02){$1$}
\uput{6pt}[d](-1.68,-0.02){$-1$}
\uput{6pt}[d](3,-0.06){$r$}
\uput[d](4,-0.02){$k$}
\uput[d](-0.8,2.6){$D_A$}
\uput[d](1.19,2.46){$f$}
\uput[d](1.7,1.23){$g$}

\end{pspicture}}}
\vspace*{-3mm}
\caption{}
\end{figure}

Denote $A^{*}:=A\backslash \left( \eta_1 \cup \eta_2 \right)$. Since $A^{*}$ is symmetric with respect to $\eta_3$, so $\tau(A^{*})$ is symmetric with respect to $\tau(\eta_3)=(-1,1)$. We reflect $D^0_{A^+}$ along $(-1,1)$ to obtain $D^0_{A^-}$, then $D^0_{A^+} \cup D^0_{A^-} \cup (-1,1)=\tau(A^{*})=:D^0_{A}$ is a fundamental domain of $A$ in $\mb{D}$. Side $\tau(\eta_1)$ is paired with its conjugate $\overline{\tau(\eta_1)}$ by a hyperbolic transformation $\tilde{f}$ fixing two points, say, $\zeta_0$ and $\bar{\zeta_0}$, and $\tau(\eta_2)$ is paired $\overline{\tau(\eta_2)}$ by a parabolic transformation $\tilde{g}$ fixing $-1$. By a M\"{o}bius transformation, under conjugacy, $\varphi: \mb{D}\rightarrow \mb{H}$ satisfying $\varphi(-1)=1$, $\varphi(1)=-1$ and $\varphi(\zeta_0)=0$, $\varphi(\bar{\zeta_0})=\infty$, $f:=\varphi \tilde{f}$, $g:=\varphi \tilde{g}$ have the form \eqref{transformations} with a parameter pair $(k,r)$, $1<r<k$, and $G=\langle f, g \rangle$ is the covering group acting on $\mb{H}$. Moreover, we let $S_1:=\{z \in \mb{H}: |z-\frac{r+1}{2r}|=\frac{r-1}{2r}\}$, $S_2:=\{z \in \mb{H}: |z|=\frac{1}{k}\}$, and let $D_A$ be the subdomain of $\mb{H}$ bounded by $S_1$, $g(S_1)$, $S_2$, $f(S_2)$ and three Euclidean line segments $(-k, -k^{-1})$, $(k^{-1}, r^{-1})$, $(r,k)$. Then $D_A^0$ is mapped onto the fundamental domain $D_A$ of $G$ in $\mb{H}$, shown in Figure \ref{fund domain in upper}. 
By the symmetry of $A$ we know $\gamma_1$ is orthogonal to $\eta_1$, $\gamma_1$ is a lift of the hyperbolic line with two endpoints $\zeta_0$ and $\bar{\zeta_0}$ in $\mb{D}$ under $\pi$, which means that $\gamma_1$ is corresponding to $\tilde{f}$ acting on $\mb{D}$, thus $f$ on $\mb{H}$.
\hfill $\Box$\vspace*{2pt}

The concept of the extremal length can be established as follows. Let $\Omega \subseteq \mb{C}$ and $\Gamma$ be a collection of finite unions of curves in $\Omega$. All of the metrics which are conformal with respect to the Euclidean metric can be defined in terms of a density $\varrho(z)|dz|$ where $\varrho(z)$ is a non-negative Borel measurable function on $\Omega$. Then the length of $\gamma \in \Gamma$ and the area of $\Omega$ with respect to $\varrho(z)$ are given by
$$L(\gamma, \varrho)=\int_{\gamma} \varrho(z)|dz|, \quad A(\Omega, \varrho)=\int_{\Omega} \varrho(z)^2 dxdy$$
with $z=x+iy$. These two quantities do not change under conformal mappings. We let $L(\Gamma, \varrho)=\inf_{\gamma \in \Gamma} L(\gamma, \varrho)$. The \emph{extremal length} of $\Gamma$ in $\Omega$ is defined by
\be \label{def of ex leng}
\lambda(\Gamma)=\sup_{\varrho} \frac{L(\Gamma, \varrho)^2}{A(\Omega, \varrho)},\ee
where the supremum is taken over all conformal densities such that $0<A(\Omega, \varrho)<\infty$. The extremal length is a conformal invariant and does not change when $\varrho$ is multiplied by a constant (see \cite[7.6.2]{Keenbook}).

Let
\ben
K(r)=\int^1_0 \frac{dx}{\sqrt{(1-x^2)(1-r^2x^2)}}
\een
with $0<r<1$ be Legendre's complete elliptic integrals of the first kind, and
\ben
\sn(u,r)=\tau \ \ \mathrm{where} \ \ u=\int^{\tau}_0 \frac{dx}{\sqrt{(1-x^2)(1-r^2x^2)}}
\een
be the Jacobian elliptic sine function. Function $\sn(u,r)$ is a bijection from $[-K(r), K(r)]$ onto $[-1,1]$ (see \cite{GDA1}). Two other functions are then defined by
\be \label{cn and dn}
\cn (u,r)=\sqrt{1-\sn^2(u,r)}, \quad \dn (u,r)=\sqrt{1-r^2\sn^2(u,r)}.
\ee
The three functions $\sn(u,r)$, $\cn (u,r)$ and $\dn (u,r)$ are called \emph{Jacobian elliptic functions} (see \cite{Byrdhandbook} for the fundamental relations and addition formulas of them). The parameter $r \in (0,1)$ is called the \emph{modulus} and the \emph{complementary modulus} of $r$ is $r'=\sqrt{1-r^2}$. If the modulus is fixed we can write Jacobian elliptic functions as $\sn \,u$, $\cn \,u$ and $\dn \,u$ for short. In the rest of this article we use $K(r)$ to refer to the Legendre's complete elliptic integral of the first kind and denote $K'(r)=K(r')=K(\sqrt{1-r^2})$. We define the \emph{normalized quotient function}
$$\mu(r)=\frac{\pi}{2}\frac{K'(r)}{K(r)}$$ for $0<r<1$, then $\mu(r)$ is a strictly decreasing homeomorphism of the interval $(0,1)$ onto $(0,\infty)$ with limit values $\mu(0+)=\infty$, $\mu(1-)=0$ (see \cite{GDA1}).

\section{Hyperbolic structure of $A$}
To prove Theorem \ref{main theorem} we need the following two theorems about the connections between $(\theta_1,\theta_2)$ and $(k,r)$, $(l_1,l_2)$ and $(k,r)$.

\begin{theorem} \label{collars}
\textsl{Suppose that $\theta_1$ and $\theta_2$ are the angular widths of the maximal peripheral collars about $\gamma_1$, $\gamma_2$. Then we have
\be\label{theta}
\cos \theta_1=\frac{r-1}{r+1},\ \ \ \cos \theta_2=\frac{t-1}{t+1}=\frac{2r(r+1)-2\delta}{\delta(r+1)-(r+1)^2},
\ee
where
\be\label{t and delta}
t=\frac{(r-1)(r+1+\delta)}{(r+3)\delta-(r+1)(3r+1)},\ \ \ \delta=k^2+r-\sqrt{(k^2-1)(k^2-r^2)}
\ee
with $k$ and $r$ as in (\ref{transformations}).}
\end{theorem}
\textbf{Proof.} In the fundamental domain $D_A$ shown in Figure \ref{fund domain in upper}, we take $S_1=\{z \in \mb{H}: |z-\frac{r+1}{2r}|=\frac{r-1}{2r}\}$. It is clear that the common tangent Euclidean line of $S_1$ and $g(S_1)$ is going through the origin. Denote the segment of the common tangent line in $\mb{H}$ by $L$ and the tangent point where $L$ and $S_1$ are tangent by $P_1$, the tangent point where $L$ and $g(S_1)$ are tangent by $P_2$. We assert that the maximal peripheral collar about the axis of $f$ is the region $C$ bounded by the non-positive real axis and $L$ in $\mb{H}$. It can be seen as follows. From the foundation of elementary geometry we know that the center of the Euclidean circle $\mathcal{C}_0$ passing through $P_1$, $P_2$ and $1$ is in $L$, and $\mathcal{C}_0$ is tangent to the real axis at $1$. Then $\mathcal{C}_0$ is a horocycle orthogonal to $S_1$ and $g(S_1)$, and it is invariant under $g$. Hence $P_2=g(P_1)$, which means $g(S_1)$ is tangent to $g(L)$ at $P_2$ if we note that $g$ is one-to-one on $\mb{H}$, and then $L$ is tangent to $g(L)$ at $P_2$, where $g(L)$ is the hypercircle in the exterior of $D_A$ with two endpoints in the interval $(1,r)$. This verifies the maximality of the collar $C$. To describe $C$ we consider the angular width $\theta$ between $L$ and the positive imaginary axis.

\begin{figure}[htbp]
\centering
 \subfigure[\label{theta1}]{\scalebox{0.65}[0.65]{
\begin{pspicture}(-4.8,-0.64)(4.8,4.9)
\psarc[fillstyle=solid,fillcolor=hylightgray,linecolor=white](0,0){4.4}{0}{180}
\psarc[fillstyle=solid,fillcolor=white,linecolor=white](0,0){1.2}{0}{180}
\psarc[fillstyle=solid,fillcolor=white,linecolor=white](1.54,0){0.34}{0}{180}
\psarc[fillstyle=solid,fillcolor=white,linecolor=white](3.6,0){0.8}{0}{180}

\psline(-4.8,0)(4.8,0)
\psline(0,0)(0,4.8)
\psline[linestyle=dashed](0,0)(4.4,1.02)
\psarc(0,0){0.24}{10}{90}
\psarc(0,0){4.4}{0}{180}
\psarc(0,0){1.2}{0}{180}
\psarc(1.54,0){0.34}{0}{180}
\psarc(3.6,0){0.8}{0}{180}
\psline{->}(0.72,1.15)(2.24,3.58)
\psarcn{->}(2.64,-0.56){1.44}{133}{70}

\uput[d](-4.4,0){$-k^2$}
\uput{5pt}[d](-1.2,0){$-1$}
\uput{5pt}[d](1.2,0){$1$}
\uput{6pt}[d](1.84,-0.04){$r$}
\uput[d](2.8,0){${k^2}/{r}$}
\uput[d](0,0){$0$}
\uput[d](4.4,0){$k^2$}
\uput[d](1.3,2.8){$f$}
\uput[d](2.48,1.48){$fg^{-1}$}
\uput[ur](0.08,0.08){$\theta_1$}
\end{pspicture}}}\ \ \hspace{3mm}
  \subfigure[\label{theta2}]{\scalebox{0.63}[0.63]{
\begin{pspicture}(-5.6,-0.64)(5.6,5.6)
\psarc[fillstyle=solid,fillcolor=hylightgray,linecolor=white](0,0){5.53}{0}{180}
\psarc[fillstyle=solid,fillcolor=white,linecolor=white](0,0){1.29}{0}{180}
\psarc[fillstyle=solid,fillcolor=white,linecolor=white](0.99,0){0.3}{0}{180}
\psarc[fillstyle=solid,fillcolor=white,linecolor=white](1.86,0){0.57}{0}{180}
\psarc[fillstyle=solid,fillcolor=white,linecolor=white](4.21,0){1.32}{0}{180}

\psline(-5.84,0)(5.84,0)
\psarc(0,0){5.53}{0}{180}
\psarc(4.21,0){1.32}{0}{180}
\psarc{->}(1.47,0.24){0.48}{70}{160}
\psarcn{->}(3.1,0.04){1.4}{150}{70}
\psarc(0,0){1.29}{0}{180}
\psarc(0.99,0){0.3}{0}{180}
\psarc(1.86,0){0.57}{0}{180}
\psline[linestyle=dashed](0,0)(5.55,1.85)
\psline(0,0)(0,5.94)

\uput[d](1.23,-0.05){$1$}
\uput[d](-1.23,0){$T(r)$}
\uput[d](0,-0.05){$0$}
\uput[d](0.86,6.2){$T(L_{fg^{-1}})$}
\uput[d](1.32,1.18){$g$}
\uput[d](2.48,1.9){$f$}
\uput[d](0.27,0.7){$\theta_2$}
\uput[d](2.2,0){$T(-1)$}
\uput[d](3.4,0){$T(-k^2)$}
\uput[d](5.53,0){$T(k^2)$}
\uput[d](-5.53,0){$T(\frac{k^2}{r})$}
\end{pspicture}}}
\vspace*{-3mm}
\caption{}
\end{figure}

For $f$ and $g$ given by \eqref{transformations}, we have
\be \label{fg inverse}
fg^{-1}=\frac{\pm 1}{r-1}\left(
\begin{array}{cc}
2k& -k(r+1)  \\
(r+1)k^{-1}& -2rk^{-1}
\end{array}
\right),
\ee
which is a hyperbolic transformation. Since $\gamma_1$, $\gamma_2$ are corresponding to $f$, $fg^{-1}$, $\theta_1$ and $\theta_2$ are angular widths of the maximal peripheral collars about the axes of $f$, $fg^{-1}$, respectively. $\theta_1$ is shown in Figure \ref{theta1}, thus $\cos\ \theta_1=\frac{r-1}{r+1}.$ Now we consider $\theta_2$, which is shown as in Figure \ref{circle collar}, where $x_1$ and $x_2$ are the fixed points of $fg^{-1}$. \begin{figure}[htbp]
\centering \scalebox{0.63}[0.63]{
 \begin{pspicture}(-9.5,-0.09)(6,7.3)
\psarc[fillstyle=solid,fillcolor=hylightgray,linecolor=white](-1.8,0){7.3}{0}{180}
\psarc[fillstyle=solid,fillcolor=white,linecolor=white](-1.8,0){3.3}{0}{180}
\psarc[fillstyle=solid,fillcolor=white,linecolor=white](2.1,0){0.6}{0}{180}
\psarc[fillstyle=solid,fillcolor=white,linecolor=white](4.5,0){1}{0}{180}
\psline(-9.5,0)(6,0)
\psarc(-1.8,0){7.3}{0}{180}
\psarc(-1.8,0){3.3}{0}{180}
\psarc(2.1,0){0.6}{0}{180}
\psarc(1.78,0){0.28}{0}{180}
\psarc(4.5,0){1}{0}{180}
\psarc(3.2,0){0.9}{0}{180}
\uput[d](-9.1,0){$-k^2$}
\uput{7pt}[d](-5.1,0.1){$-1$}
\uput{7pt}[d](1.5,0.1){$1$}
\uput{8pt}[d](2.73,0.1){$r$}
\uput[d](3.57,0.1){$\frac{k^2}{r}$}
\uput[d](5.5,0.1){$k^2$}

\psarc[linestyle=dashed](3.19,1.8){2}{-62}{240}
\psline(2.52,0.62)(2.66,0.49)(2.55,0.36)
\psline(1.7,0.2)(1.9,-0.6)
\psline(3.64,0.78)(3.55,0.66)(3.68,0.55)

\uput{8pt}[d](2.3,0.1){$x_1$}
\uput{8pt}[d](4.3,0.1){$x_2$}
\uput[d](2.17,0.7){$\theta_2$}
\uput[d](4.32,0.86){$\theta_2$}
\uput[d](3.18,1.56){$L_{fg^{-1}}$}
\uput[d](-3.8,3.3){$S$}
\uput[d](1.96,-0.5){$g(S)$}
\uput[d](-2.96,5.5){$D'_A$}
\uput[d](4.32,0.86){$\theta_2$}
\end{pspicture}}
\caption{\label{circle collar}}
\end{figure}By some M\"{o}bius transformation $T$ with $T(x_1)=0$, $T(x_2)=\infty$, $D'_A$ in Figure \ref{circle collar} can be mapped onto the fundamental domain shown in Figure \ref{theta2} with the vertices $T(-k^2)$, $T(-1)$, $T(1)$, $T(r)$, $T(x_2)$, $T(k^2)$ and the angle $\theta_2$. To identify $T$, without loss of generality, we may assume that $T(1)=1$. Then after some computation, we have
$$x_1=\frac{\delta}{r+1} \in \left(\frac{r}{2}+\frac 1 2, \ r\right), \quad x_2=\frac{k^2}{x_1}.$$
So
$$T(z)=\frac{(1-r)z-x_1(1-r)}{(r+1-2x_1)z+x_1(r+1)-2r}.$$
Thus
$$T(-1)=\frac{(r-1)(r+1+\delta)}{(r+3)\delta-(r+1)(3r+1)}:=t.$$
We have $t>1$ provided that $2r<\frac 1 2 (r+1)^2 <\delta < r(r+1)$. That is shown in Figure \ref{theta2} and then $\displaystyle \cos \ \theta_2=\frac{t-1}{t+1}. $ \vspace*{1mm}\hfill $\Box$

{\setlength{\parindent}{0pt}
\textbf{Remark 1.}} Theorem  \ref{collars} is an improvement of the collar lemma. The maximal collar defined by \eqref{kl collar} with equality is smaller than the collar decided by $\theta_1$ and $\theta_2$ in \eqref{theta}. We denote the angular widths of the collars defined by \eqref{kl collar} about the axes of $f$ and $fg^{-1}$ by $\theta'_1$ and $\theta'_2$. Then $\cos \theta'_1=\frac{k^2-1}{k^2+1}$ by \eqref{kl collar}. Since the function $\frac{x-1}{x+1}$ is monotonically increasing for $x \in [1, \infty)$, and $1 <r <k<k^2$, then $\theta'_1 < \theta_1$ from the expression of $\cos \theta_1$ in \eqref{theta}. From the symmetry given by \eqref{main quantity}, or by pure computation, we have $\theta'_2 < \theta_2$ for $fg^{-1}$.

\begin{lemma}\label{connection l and kr}
\textsl{In the punctured annulus $A$, for the lengths $l_1$, $l_2$ of $\gamma_1$, $\gamma_2$, and $(k,r)$ as in \eqref{transformations}, we have
\be \label{l1 and l2}
2\cosh \frac{l_1}{2}=k+\frac{1}{k},\ \ \ 2\cosh \frac{l_2}{2}=\frac{2}{r-1}\left(k-\frac{r}{k}\right).
\ee}
\end{lemma}
\textbf{Proof.} Since $\gamma_1$, $\gamma_2$ are corresponding to $f$, $fg^{-1}$, respectively, $l_1$ and $l_2$ are the translation lengths of $f$ and $fg^{-1}$. Thus \eqref{l1 and l2} is obvious from \eqref{fg inverse}. \hfill $\Box$ \vspace*{2mm}

{\setlength{\parindent}{0pt}
\textbf{Proof of Theorem \ref{main theorem}.}} At first we solve that $$r=\frac{k\cosh \frac{l_2}{2} +k^2}{k\cosh \frac{l_2}{2} +1}$$ from the second equation of \eqref{l1 and l2}, then by \eqref{theta} and the first equation of \eqref{l1 and l2},
$$\cos \theta_1= \frac{r-1}{r+1}=\frac{k-\frac{1}{k}}{2\cosh \frac{l_2}{2}+k+\frac{1}{k}}=\frac{\sinh \frac{l_1}{2}}{\cosh \frac{l_1}{2}+\cosh\frac{l_2}{2}}$$
as required. The expression for $\cos \theta_2$ in \eqref{main quantity} can be obtained by symmetry. \hfill $\Box$ \vspace*{1mm}

{\setlength{\parindent}{0pt}
\textbf{Remark 2.}} We compare $\theta_1$ and $\theta_2$, $l_1$ and $l_2$. By \eqref{main quantity} and
$$t-r=\frac{(r+1)^2(3r-1-\delta)}{(r+3)\delta-(r+1)(3r+1)}, \quad -(r-1)^2<3r-1-\delta<\frac 1 2 (r-1)(r-3),$$
we have the following inequalities. When $1<r<3$,
$$ t>r,\  \theta_1>\theta_2\ \mathrm{and}\ l_1<l_2,\qquad \qquad \qquad \qquad \mbox{if\ }\ \sqrt{\frac{3r-1}{3-r}}<k,$$
$$ t=r,\  \theta_1=\theta_2\ \mathrm{and}\ l_1=l_2,\qquad \qquad \qquad \qquad \mbox{if\ }\ \sqrt{\frac{3r-1}{3-r}}=k,$$
$$ t<r,\  \theta_1<\theta_2\ \mathrm{and}\ l_1>l_2,\qquad \qquad \qquad \qquad \mbox{if\ }\ k<\sqrt{\frac{3r-1}{3-r}};$$
when $r \geq 3$,\ $t<r$, $\theta_1<\theta_2$ and $\ l_1>l_2$. Then the following corollary is obtained provided that $a=1$ is the midpoint between $1/R$ and $R$ in the hyperbolic metric on $A$.
\begin{corollary} \label{coro}
\textsl{In the punctured annulus $A=\{z:1/R<|z|<R\} \backslash \{1\}$, the covering group $G$ of $A$ is generated by
\ben
f(z)= \frac{3r-1}{3-r}z,\ \ \ g(z)= \frac{2rz-(r+1)}{(r+1)z-2},
\een
where $1<r<3$ and $r$ is related to $R$ in some unknown way.}
\end{corollary}

{\setlength{\parindent}{0pt}
\textbf{Remark 3.}} Here we provide another way to prove Theorem \ref{main theorem} by the use of the pants decomposition for a hyperbolic manifold with a cusp. It is different from the description presented in this article, but it is meaningful for the discussion of pants partitions and spectral questions. We only show the rough idea here without rigorous proof. For more details on pairs of pants and collars (see \cite{buserbook}). Note that $\gamma_1$ and $\gamma_2$ are orthogonal to the real axis, at points, say, $s_1$ and $s'_1$, $s_2$ and $s'_2$, respectively. Let $\widetilde{C}_{\theta_1}(\gamma_1)$ be the maximal peripheral collar about geodesic $\gamma_1$ of angular width $\theta_1$. Then, by symmetry, the inner boundary $l$ of $\widetilde{C}_{\theta_1}(\gamma_1)$ is orthogonal to segments $(-R, -1/R)$, $(a,R)$, and tangent to itself at a point between $1/R$ and $a$. The sketch of $l$, $\gamma_1$ and $\gamma_2$ is shown in Figure \ref{ell shape}. \begin{figure}[htbp]
\centering \scalebox{0.71}[0.71]{
 \begin{pspicture}(-3,-3)(3,2.8)
\pscircle[fillstyle=solid,fillcolor=hylightgray,linecolor=white](0,0){3}
\pscircle[fillstyle=solid,fillcolor=white,linecolor=white](0,0){0.5}
\pscircle(0,0){3}
\pscircle(0,0){0.5}
\psdot[dotsize=3pt](1.6,0)
\psdot[dotsize=2pt](-2.6,0)
\psdot[dotsize=2pt](-1.7,0)
\psdot[dotsize=2pt](0.6,0)
\psdot[dotsize=2pt](2.4,0)
\psline[linewidth=0.4pt](-3,0)(-0.5,0)
\psline[linewidth=0.4pt](3,0)(0.5,0)
\psccurve[linestyle=dashed](2.4,0)(1.3,1.9)(0.3,2.2)(-0.7,1.9)(-1.7,1.3)(-2.6,0)(-1.7,-1.3)(-0.7,-1.9)(0.3,-2.2)(1.3,-1.9)(2.4,0)
\psccurve[linestyle=dashed](0.6,0)(0.3,0.6)(-0.2,0.7)(-0.7,0.65)(-1.2,0.5)(-1.7,0)(-1.2,-0.5)(-0.7,-0.65)(-0.2,-0.7)(0.3,-0.6)(0.6,0)
\psccurve[linewidth=0.4pt](2.1,0)(1.8,0.4)(1.5,0.4)(1.2,0.01)(0.8,0.8)(0.4,1)(-0.2,0.9)(-0.8,0.4)(-0.9,0)(-0.8,-0.4)(-0.2,-0.9)(0.4,-1)(0.8,-0.8)(1.2,-0.01)(1.5,-0.4)(1.8,-0.4)(2.1,0)
\uput[d](1.6,0){$a$}
\uput[d](1.3,2.6){$\gamma_1$}
\uput[u](1,0.5){$l$}
\uput[u](-1.1,0.5){$\gamma_2$}
\uput[d](2.65,0.08){$s'_1$}
\uput[d](-2.78,0.02){$s_1$}
\uput[d](0.76,0.08){$s'_2$}
\uput[d](-1.8,0.02){$s_2$}
\end{pspicture}}
\caption{\label{ell shape}}
\end{figure} The punctured domain bounded by $\gamma_1$ and $\gamma_2$ with puncture $a$ is conformally equivalent to a pair of pants with a cusp, called a Y-piece with a cusp (see \cite{buserbook}), shown in Figure \ref{v piece}. We use the same notations as in $A$. Then, in the Y-piece, $l$ is a curve orthogonal to geodesic ray $(a, s'_1)$ and tangent to $(a, s'_2)$, and also orthogonal to the geodesic line containing $(s_1, s_2)$. It is known that the Y-piece in Figure \ref{v piece} with a cusp can be decomposed into two isometric pentagons with four right angles, one of which is shown in Figure \ref{pentagon}, where $l$ is corresponding to a hypercycle, still denoted by $l$, tangent to $(a, s_2)$. It is easy to see that, in this pentagon, the angle between $l$ and the geodesic containing side $\frac 1 2 \gamma_1$ is the angular width $\theta_1$ of $\widetilde{C}_{\theta_1}(\gamma_1)$. So that we can obtain \eqref{main quantity} for $\theta_1$ only from properties of a right-angled pentagon, and then for $\theta_2$ by symmetry. Moreover, the hyperbolic length of segment $(s_1, s_2)$ in Figure \ref{pentagon} can be obtained,
$$\sinh  l(s_1, s_2) =\frac{\cosh \frac{l_1}{2}+\cosh \frac{l_2}{2}}{\sinh \frac{l_1}{2}\;\!\sinh \frac{l_2}{2}},$$
it is the distance between  $\gamma_1$ and $\gamma_2$ in Figure \ref{ell shape}.

\begin{figure}[htbp]
 \centering
 \subfigure[\label{v piece}]{\scalebox{0.96}[0.96]{
\begin{pspicture}(-2.5,-4.5)(2.5,0.5)

\pscurve(-0.03,0.07)(-0.03,0)(-0.21,-1)(-0.8,-2.1)(-1.5,-3)
\pscurve(0.03,0.07)(0.03,0)(0.21,-1)(0.8,-2.1)(1.5,-3)
\psccurve(-1.5,-3)(-1,-3)(-0.7,-3.3)(-0.7,-3.7)(-1,-3.78)(-1.5,-3.4)
\psccurve(1.5,-3)(1,-3)(0.7,-3.3)(0.7,-3.7)(1,-3.78)(1.5,-3.4)
\pscurve(-0.64,-3.6)(-0.5,-3.2)(-0.3,-2.9)(0,-2.76)(0.3,-2.9)(0.5,-3.2)(0.65,-3.6)
\pscurve[linewidth=0.4pt](0.97,-3.03)(0.98,-2.9)(0.88,-2.4)(0.64,-1.9)(0.36,-1.62)(0,-1.46)(-0.2,-1.42)(-0.3,-1.41)(-0.4,-1.43)
\pscurve[linestyle=dashed,linewidth=0.4pt](-0.4,-1.46)(-0.3,-1.61)(-0.2,-1.68)(0,-1.8)(0.4,-1.86)(0.62,-1.9)(0.7,-2.4)(0.74,-2.9)(0.84,-3.23)(0.98,-3.5)(1.12,-3.7)

\psline[linewidth=0.4pt](-0.3,-1.24)(-0.1,-1.24)(-0.18,-1.43)
\psline[linewidth=0.4pt](-1.4,-3.6)(-1.6,-3.96)
\psline[linewidth=0.4pt](1.4,-3.6)(1.6,-3.96)

\psdot[dotsize=2pt](-1.5,-3)
\psdot[dotsize=2pt](-0.7,-3.7)
\psdot[dotsize=2pt](0.7,-3.7)
\psdot[dotsize=2pt](1.5,-3)

\uput[d](0,0.56){$a$}
\uput[d](-1.6,-3.86){$\gamma_1$}
\uput[d](1.62,-3.86){$\gamma_2$}
\uput[d](0.4,-1.76){$l$}

\uput[d](-1.72,-2.7){$s'_1$}
\uput[d](1.76,-2.7){$s'_2$}
\uput[d](-0.54,-3.56){$s_1$}
\uput[d](0.54,-3.5){$s_2$}

\end{pspicture}}  }\ \ \hspace{2cm}
\subfigure[\label{pentagon}]{\scalebox{0.76}[0.76]{
\begin{pspicture}(-3,-3)(3,2.6)

\pscircle(0,0){3}
\psdot[dotsize=3pt](0,2)
\psdot[dotsize=3pt](0,0)
\psdot[dotsize=3pt](1.25,0)
\psdot[dotsize=3pt](1.66,1.52)
\psdot[dotsize=3pt](1.15,2.75)
\psline(0,0)(0,2)
\psline[linestyle=dashed](0,3)(0,2)
\psline[linestyle=dashed](0,0)(0,-3)
\psline(0,0)(1.25,0)
\psarc(0,3.25){1.25}{-90}{-24}
\psarc(2.2,2.46){1.1}{-197}{-119}
\psarc(4.23,0){2.98}{150}{180}
\psarc[linestyle=dashed](-1.58,0){3.39}{-62}{21}
\psarc[linestyle=dashed](-1.58,0){3.39}{45}{62}
\psarc[linewidth=0.4pt](-1.58,0){3.39}{22}{44}

\psarc[linewidth=0.4pt](0.06,2.96){0.24}{-106}{-36}

\psline[linewidth=0.3pt](0,0.2)(0.2,0.2)(0.2,0)
\psline[linewidth=0.3pt](1.05,0)(1.05,0.2)(1.26,0.2)
\psline[linewidth=0.3pt](0,1.8)(0.2,1.8)(0.2,2.01)
\psline[linewidth=0.3pt](1.52,1.62)(1.41,1.44)(1.59,1.32)
\psline[linewidth=0.3pt](0.74,2.22)(0.86,2.09)(0.99,2.22)

\uput[d](1,1.2){$\frac 1 2 \gamma_2$}
\uput[u](1,1.56){$l$}
\uput[u](-0.31,0.5){$\frac 1 2 \gamma_1$}
\uput[d](-0.2,2.4){$s'_1$}
\uput[d](-0.2,0.1){$s_1$}
\uput[d](1.84,1.76){$s'_2$}
\uput[d](1.25,0.02){$s_2$}
\uput[d](1.16,3.24){$a$}
\uput[d](0.36,2.8){$\theta_1$}

\end{pspicture}}}
\vspace*{-5mm}
\caption{}
\end{figure}

\section{Complex structure of $A$}

For the proof of Theorem \ref{extremal length theorem} we need the following three lemmas.
\begin{lemma}\textrm{$\mathrm{(see\ 6.26,\ 6.27\ in\ [1])}$} \label{halfannulua onto halfdisk}
\textsl{ For $0<q<1$ let $\K:=K(q)$,\ $\K':=K'(q)$ and select $b=\exp(-\pi \K'/(4\K))$. Then the conformal mapping $\omega$ defined by $$\omega(z)=\sq{q}\,\sn\left(\frac{2i\K}{\pi}\llog {\frac{z}{b}} +\K, q \right)$$
is unique up to rotations and takes the annulus $b<|z|<1$ onto the unit disk $|\omega(z)|<1$ minus the slit $[-\sq q, \sq q]$, where $\llog$ is the principal branch of the logarithm.}
\end{lemma}

\begin{lemma} \label{unsymmetric onto symmetric}
\textsl{For given numbers $q$ and $\alpha$, $0<\sq q<\alpha<1$, the unique M\"{o}bius transformation $\sigma$ which preserves $\mb{D}$ and satisfies $\sigma(-1)=-1$, $\sigma(1)=1$, $\sigma(-\sq q)=0$ is given by
is
$$\sigma(z)=\frac{z+\sq q}{\sq q z+1}.$$}
\end{lemma}
\textbf{Proof.} Note that the correspondences between three distinct points and their images on $\hat{\mb{C}}$ decide a M\"{o}bius transformation and a M\"{o}bius transformation preserving the unit disk has the form $\pm \left(
\begin{array}{cc}
b_1& \overline{b_2}  \\
b_2& \overline{b_1}
\end{array}
\right)\in{\textrm{PSL}(2,\mb{C})}$, then the conclusion is obvious. \hfill $\Box$ \vspace{3mm}

The actions of $\omega$ and $\sigma$ are shown in Figure \ref{VAV's fig}.

\begin{figure}[htbp]
{\scalebox{0.7}[0.7]{
\begin{pspicture}(-7.6,-3)(18,3)

\psline{->}(-1.3,0.3)(-0.1,0.3)
\uput[d](-0.75,1.3){\LARGE{$w$}}
\psline{->}(6.14,0.3)(7.44,0.3)
\uput[d](6.79,1.3){\LARGE{$\sigma$}}

\pscircle[fillstyle=solid,fillcolor=hylightgray,linecolor=white](-4.7,0){2.5}
\pscircle[fillstyle=solid,fillcolor=white](-4.7,0){0.6}
\pscircle(-4.7,0){2.5}
\pscircle(-4.7,0){0.6}
\psdot(-7.2,0)
\psdot(-2.2,0)
\psdot(-5.3,0)
\psdot(-4.1,0)
\uput[d](-7.46,0){$-1$}
\uput[d](-5.56,0){$-b$}
\uput[d](-4,0){$b$}
\uput[d](-2.06,0){$1$}

\pscircle[fillstyle=solid,fillcolor=hylightgray,linecolor=white](3,-0){2.5}
\pscircle(3,0){2.5}
\psline(2.4,0)(3.6,0)
\psdot(0.5,0)
\psdot(5.5,0)
\psdot(4.3,0)

\psline(2.24,0)(3.76,0)
\psline(2.24,0)(3.76,0)
\uput[d](0.2,0){$-1$}
\uput[d](5.64,0){$1$}
\uput[d](4.3,-0.06){$\alpha$}
\uput[d](2.3,0){$-\sq q$}
\uput[d](3.7,-0.06){$\sq q$}

\pscircle[fillstyle=solid,fillcolor=hylightgray,linecolor=white](10.7,0){2.5}
\pscircle(10.7,0){2.5}
\psdot(10.7,0)
\psdot(12.15,0)
\psdot(8.2,0)
\psdot(13.2,0)
\psdot(11.75,0)
\psline(10.7,0)(11.75,0)
\psline(10.7,0)(11.75,0)
\psline(10.7,0)(11.75,0)
\uput[d](7.9,0){$-1$}
\uput[d](13.36,0){$1$}
\uput[d](10.7,0){$0$}
\uput[d](11.6,-0.04){$\sigma(\sq q)$}
\uput[d](12.4,0.58){$\sigma(\alpha)$}

\end{pspicture}}}
\caption{\label{VAV's fig}}
\end{figure}

\begin{lemma} \label{ohtsuka theorem}
\textsl{Let $\widetilde{\mc{C}}$ be the family of loops in $\mathbb{D}$ separating $0$ and $x$ from the unit circle $\partial \mathbb{D}$, $0<x<1$, and ${\mc{C}}$ be the family of loops in $\mathbb{D}$ separating the slit $(0, x)$ from $\partial \mathbb{D}$. Then the extremal lengths of $\widetilde{\mc{C}}$ and ${\mc{C}}$ are
\be \label{lambdaC}
\lambda(\widetilde{\mc{C}})=\lambda(\mc{C})=\frac{2\pi}{\mu(x)}.
\ee}
\end{lemma}

Lemma \ref{ohtsuka theorem} can be found in \cite[2.1]{Lehto Virtanen} for the extremal length of family $\mc{C}$ as the reciprocal of the modulus of Gr\"{o}tzsch's extremal domain. On the twice-punctured unit disk $\mb{D} \backslash \{\pm r\}$ Ohtsuka \cite[Theorem 2.56]{OHT} gave the extremal length of the curve family which separates punctures $\{-r, r\}$ from the unit circle. If we note that there exists a M\"{o}bius transformation which maps Ohtsuka's domain $\mb{D} \backslash \{\pm r\}$ to the domain $\mb{D} \backslash \{0,x\}$ and $\mu(\frac{2 r}{r^2+1})=\frac{1}{2}\mu(r^2)$, we can obtain the extremal length of $\widetilde{\mc{C}}$ in \eqref{lambdaC}.\\[2mm]
\textbf{Proof of Theorem \ref{extremal length theorem}.} At first, we take the mappings $$\varphi_1(z)=\frac{z}{R}, \quad \varphi_2(z)=\frac{1}{Rz},$$
then $\varphi_1(A),\ \varphi_2(A)$ are still punctured annuli
$$\varphi_1(A)=\left\{ z:\frac{1}{R^2}<|z|<1\right\} \setminus \left\{\frac a R \right\}, \quad \varphi_2(A)=\left\{z:\frac{1}{R^2}<|z|<1\right\}\setminus \left\{\frac{1}{a R}\right\}.$$
Next let $b=1/R^2$ and define  \ben
\omega (z)=\sq{q}\,\sn\left(\frac{2i\K}{\pi}\llog {R^2z} +\K, q \right),
\een
where $q$ satisfies $\mu(q)=4\log R$. Then from Lemma \ref{halfannulua onto halfdisk}, $\omega(z)$ maps $\varphi_1(A)$ and $\varphi_2(A)$ onto two punctured slit disks $\mb{D} \backslash \left([-\sq q, \sq q]\cup \{\alpha_1\}\right)$ and $\mb{D} \backslash \left([-\sq q, \sq q]\cup \{\alpha_2\}\right)$, respectively, where
$$\alpha_1=\omega \left(\frac{a}{R}\right)=\sq{q}\,\sn\left(iu_1+\K, q \right),\ \ \  \alpha_2=\omega \left(\frac{1}{aR}\right)=\sq{q}\,\sn\left(iu_2 +\K, q \right)$$
with $u_1$ and $u_2$ given by \eqref{u1 u2}. We note that
$$\sn\left(iu_j+\K, q \right)=\frac{\cn(iu_j, q)}{\dn(iu_j,q)},$$
and
$$\cn(iu_j, q)=\frac{1}{\cn(u_j,{q}')},\ \ \dn(iu_j, q)=\frac{\dn(u_j, {q}')}{\cn(u_j,{q}')}$$
for $j=1,2$ (see $(120.02)$, $(122.03)$ and $(125.02)$ in \cite{Byrdhandbook}). Then $\alpha_1=\frac{\sq q}{\dn\,u_1}$, $\alpha_2=\frac{\sq q}{\dn\,u_2}$ with the modulus of $\dn$ being ${q}'=\sqrt{1-q^2}$. So for $j=1,2$,
$${p_j}:=\sigma(\alpha_j)=\frac{\alpha_j+\sq q}{\sq q\alpha_j+1}=\frac{\sq q(\dn\,u_j+1)}{q+\dn\,u_j}.$$ For $\widetilde{\mc{C}}$ and ${\mc{C}}$ defined in Theorem \ref{extremal length theorem}, since $\widetilde{\mc{C}}\supset \mc{C}_j \supset \mc{C}$, then $\lambda(\widetilde{\mc{C}})\leq \lambda(\mc{C}_j) \leq \lambda({\mc{C}})$, $j=1,2$. Thus \eqref{extremal length results} is obtained from \eqref{lambdaC}. \hfill $\Box$
\begin{corollary} \label{alpha=1}
\textsl{In the same assumption as Theorem \ref{extremal length theorem}, when the puncture $a=1$,
$$\lambda_1=\lambda_2=\frac{2\pi}{\mu(p)},$$
where
\be \label{special p}
p=\frac{\sq q\left(\sqrt{\sq q(\sqrt{1+q}+1)}+\sqrt{\sqrt{1+q}+\sq q}\right)}{q\sqrt{\sqrt{1+q}+\sq q}+\sqrt{\sq q(\sqrt{1+q}+1)}}.
\ee}
\end{corollary}
\textbf{Proof.} When the puncture $a=1$, we note that $u_1=u_2=\frac{\K'}{4}$, and by (122.10) in \cite{Byrdhandbook}, $$\sn\frac{\K'}{2}=\frac{1}{\sqrt{1+q}},\ \ \cn\frac{\K'}{2}=\frac{\sq q}{\sqrt{1+q}},\ \ \dn\frac{\K'}{2}=\sq q,$$
where all Jacobian elliptic function $\sn$, $\cn$ and $\dn$ have the modulus ${q}'$, then the half argument formulas (124.01) in \cite{Byrdhandbook} lead to \ben
\dn \frac{\K'}{4}=q^{\frac{1}{4}}(\sqrt{1+q}+1)^{\frac{1}{2}}(\sqrt{1+q}+\sq q)^{-\frac{1}{2}}.\een
Combining with \eqref{p1p2} we obtain \eqref{special p}.
\hfill $\Box$ \vspace{2mm}

We can compare the hyperbolic and extremal lengths.
\begin{theorem} \label{comparison of hyper and ex}
\textsl{With the same $\lambda_j$ as in \eqref{extremal length results}, $l_j$ as in \eqref{l1 and l2}, $\theta_j$ as in \eqref{theta}, $j=1,\,2$, we have
\be \label{comparison ex}
\displaystyle \frac{l_j}{\pi}  \leq \lambda_j \leq \frac{l_j}{\frac{\pi}{2}+\arccos \left(\frac{\sinh (l_j /2)}{\cosh(l_1 /2)+\cosh(l_2 /2)}\right)}.
\ee}
\end{theorem}
\textbf{Proof.} From \cite[(1) and (2)]{Maskit1985}, in the punctured annulus $A$ we have $(\frac{\pi}{2}+\theta_j)\lambda_j \leq l_j \leq \pi \lambda_j$, that is $\frac{l_j}{\pi} \leq \lambda_j \leq \frac{l_j}{\frac{\pi}{2}+\theta_j}$. Substituting \eqref{main quantity} leads to \eqref{comparison ex}. \hfill $\Box$

\section{Degenerating cases}
In this section we consider the degenerating cases when the puncture $a$ is approaching the inner or outer boundary of $A$ and give some connection between the parameter pairs $(R,a)$ and $(k,r)$, and some other cases when one of the boundaries is shrinking to a point or to the other boundary.

Case (i): $R$ is fixed. When $a\rightarrow R$, by the expression \eqref{u1 u2}, we know $u_1\rightarrow \frac{\K '}{2}$, $u_2\rightarrow 0$ as $a\rightarrow R$, then $\dn\,u_1 \rightarrow \sq q$, $\dn\,u_2\rightarrow 1$, and $p_1 \rightarrow 1$, $p_2 \rightarrow \frac{2\sq q}{q+1}$ (see \cite[p.\;19]{Byrdhandbook}). We note that $\mu(\frac{2\sq q}{q+1})=\frac{1}{2}\mu(q)=2 \log R$ (see \cite[5.2]{GDA1}), thus $\lambda_1 \rightarrow +\infty$, $\lambda_2 \rightarrow \frac{\pi}{\log R}$ with $\lambda_1$ and $\lambda_2$ given by \eqref{extremal length results}. From \eqref{comparison ex} we know $l_1 \rightarrow \infty$ for $\lambda_1\rightarrow \infty$, and thus $k \rightarrow +\infty$ from the first formula in \eqref{l1 and l2}. When $a\rightarrow 1/R$, it holds $\lambda_1 \rightarrow \frac{\pi}{\log R}$, $\lambda_2 \rightarrow +\infty$. Thus $l_2\rightarrow \infty$, and $r \rightarrow 1$ from the second formula in \eqref{l1 and l2}.

Case (ii): $a$ is fixed. When $R\rightarrow +\infty$, $A$ is tending to the thrice-punctured sphere. Since $\mu(q)=4\log R$, we know $q\rightarrow 0$, and $p_1 \rightarrow 0$, $p_2 \rightarrow 0$ from \eqref{p1p2}. Thus $\lambda_1 \rightarrow 0$, $\lambda_2 \rightarrow 0$, which imply that $\l_1 \rightarrow 0$, $\l_2 \rightarrow 0$ from \eqref{comparison ex}. By the first formula of \eqref{l1 and l2}, we have $k \rightarrow 1$. The second formula of \eqref{l1 and l2} is written as
$$\cosh \frac{l_2}{2}=\frac{1}{k} \frac{1}{r-1}(k^2-1-r+1)=\frac{1}{k} \left( \frac{k^2-1}{r-1}-1 \right).$$
Thus $\frac{k-1}{r-1} \rightarrow 1$, $k \rightarrow 1$ and $r \rightarrow 1$ when $R \rightarrow +\infty$. In this case from the forms of $f$ and $g$ we see that the covering group $G$ does not converges algebraically, so we need another Fuchsian group to describe the degeneration. Let $h(z):=\frac{k+1}{k-1} \frac{z-1}{z+1}$. We consider $G_1:=hGh^{-1}$. Since $h(1)=0$, $h(\frac{1}{k})=-1$, $h(k)=1$, from the representation of $f$ and $g$ in \eqref{transformations}, we have
\be
hfh^{-1}=\frac{\pm}{2k} \left(
\begin{array}{cc}
k^2+1   & (k+1)^2  \\
(k-1)^2 & k^2+1
\end{array}
\right) \rightarrow \pm \left(
\begin{array}{cc}
1  & 2  \\
0  & 1
\end{array}
\right):=f_0, \ \ \mathrm{as}\ R \rightarrow +\infty, \label{f cong}\\[1mm]
hgh^{-1}={\pm} \left(
\begin{array}{cc}
\vspace*{1mm} 1   & 0 \vspace*{1mm} \\
\vspace*{1mm} 2\frac{(k-1)(r+1)}{(k+1)(r-1)} & 1 \vspace*{1mm}
\end{array}
\right) \rightarrow \pm \left(
\begin{array}{cc}
1  & 0  \\
2  & 1
\end{array}
\right):=g_0, \ \ \mathrm{as}\ R \rightarrow +\infty. \label{g cong}
\ee
That means, when $R \rightarrow +\infty$, $h$ maps the fundamental domain $D_A^0$ shown in Figure \ref{fund domain in disk} onto the domain bounded by $\mathrm{Re}\,z=\pm 1$ and $\{|z \pm \frac{1}{2}|=\frac{1}{2}\}$. And we know that, as $R \rightarrow +\infty$, $G_1$ converges algebraically to a Fuchsian group generated by $f_0$ and $g_0$ from \eqref{f cong} and \eqref{g cong}. The group $\langle f_0,g_0\rangle$ is known as the principal congruence subgroup of $PSL(2, \mb{Z})$ of level $2$ (see e.g. \cite[p.54]{Farkaskrabook2001}). This case is the degeneration for the covering group when the two hyperbolic elements $g$, $fg^{-1}$ are both becoming parabolic ones.

Case (iii): $a=1$. When $R \rightarrow 1$, we have $q\rightarrow 1^-$ from $\mu(q)=4\log R$, thus $p_1 \rightarrow 1^-$, $p_2 \rightarrow 1^-$ and $\lambda_1 \rightarrow +\infty$, $\lambda_2 \rightarrow +\infty$. Therefore $l_1\rightarrow \infty$, $l_2\rightarrow \infty$ from \eqref{comparison ex}, and $k \rightarrow +\infty$ from the first formula in \eqref{l1 and l2}. By Corollary \ref{coro}, $k^2=\frac{3r-1}{3-r}$, we have $r\rightarrow 3$, and $g(z)$ converges to the parabolic transformation $z \mapsto \frac{3z-2}{2z-1}$. This is the degeneration for the covering group when the trace of a hyperbolic element is going to $\infty$. In the fundamental domain $D_A$ shown in Figure \ref{fund domain in upper}, the two semi-circles with end points $\pm\frac 1 k$, $\pm k$, are shrinking to two points $0$, $\infty$, respectively. Then $D_A$ is becoming the triangle $\{z \in \mb{H}: |z+\frac{1}{2r}| > \frac{1}{2r}\ \mathrm{and}\ |z-\frac{r}{2}| > \frac{r}{2}\}$. In this case, $A$ is conformally equivalent to the once-punctured unit disk.

Case (iv): $a/R$ is fixed. Now we let $R \rightarrow +\infty$ and take the once-punctured annulus model as $A_1=\{R^{-2}<|z|<1\} \backslash \{x\}$, $x=a/R$, $R^{-2}<x<1$. First we consider $\lim_{R \rightarrow +\infty} \sn (v, q')$, where $v:=\frac{2\K}{\pi}\log x$, $\K=K(q)$ and $q$ satisfies $\mu(q)=4\log R$. All the Jacobian functions $\sn$, $\cn$, $\dn$ in this section have the modulus $q'$, we will omit $q'$ in the notation of them. Since
\be \label{def of sn}
\int^{\sn\ls v}_{0} \frac{dt}{\sqrt{(1-t^2)(1-(1-q^2)t^2)}}=v,
\ee
and $q \rightarrow 0$, $\K \rightarrow \frac{\pi}{2}$ as $R \rightarrow +\infty$, taking limits of both sides of \eqref{def of sn} gives
\be \label{lim sn}
\lim_{R \rightarrow +\infty} \sn \ls v=\frac{x^2-1}{x^2+1},
\ee
and thus
\be \label{lim cn dn}
\lim_{R \rightarrow +\infty} \cn \ls v=\frac{2x}{x^2+1}, \quad \lim_{R \rightarrow +\infty} \dn \ls v=\frac{2x}{x^2+1}.
\ee
Note that
$$\dn \ls u_1=\dn \left( \frac{2\K}{\pi}\log{R^2 x}\right)=\dn\left(\frac{{\K}'}{2}+v\right)=\frac{\sqrt{q}\left(\dn\ls v-(1-q)\sn \ls v \,\cn\ls v \right)}{1-(1-q)\sn^2 v}, $$
where the last equivalence is due to \cite[122.10 and 123.01]{Byrdhandbook}. Then from \eqref{lim sn} and \eqref{lim cn dn},
\ben
\lim_{R \rightarrow +\infty} \frac{\dn \ls u_1}{\sqrt{q}}=\lim_{R \rightarrow +\infty} \frac{\dn\ls v -(1-q)\sn\ls v\, \cn\ls v}{1-(1-q)\sn^2 v}=\frac 1 x.
\een
Hence for $p_1$ as in \eqref{p1p2}, we have
\be \label{lim p1}
\lim_{R \rightarrow +\infty} p_1= \lim_{R \rightarrow +\infty} \frac{\dn\ls u_1+1}{\sqrt{q}+\dn\ls u_1 /\sqrt{q}}=\lim_{R \rightarrow +\infty} \frac{\sqrt{q}}{\dn \ls u_1}=x.
\ee
Therefore $\lambda_1 \rightarrow \frac{2\pi}{\mu(x)}$. By the second formula in \eqref{u1 u2} we see $u_2$ is fixed, so $\dn \ls u_2$ is fixed to be a finite number. Thus $p_2$ in \eqref{p1p2} is tending to $0$, and $\lambda_2 \rightarrow 0$ as $R \rightarrow +\infty$. In this case, $\lim A_1$ is conformally equivalent to the twice-punctured unit disk $\mb{D}\backslash \{0,x\}$, and $\lim_{R \rightarrow +\infty}\lambda_1$ coincides with $\lambda(\widetilde{\mc{C})}$ in \eqref{lambdaC}. To show the degeneration of $A_1$ in terms of the deformation of the covering group $G$, we let $h_1(z):=\frac{r+1}{r} \frac{z-r}{z-1}$ and consider $h_1Gh_1^{-1}$. Then $h_1(1)=\infty$, $h_1(\frac{1}{r})=\frac{(r+1)^2}{r}$, $h_1(r)=0$. And $h_1fh_1^{-1}$, $h_1gh_1^{-1}$, $h_1fg^{-1}h_1^{-1}$ satisfy
\ben
h_1fh_1^{-1}=\frac{\pm 1}{k(r^2-1)} \left(
\begin{array}{cc}
 \vspace*{1mm}  \displaystyle -(r+1)\left(k^2- r \right)   &  \displaystyle (r+1)^2\left(k^2 - 1 \right) \\[1mm]
 \displaystyle -r\left(k^2-1\right) &  \displaystyle(r+1)\left(k^2r- 1 \right)
\end{array}
\right) \rightarrow \pm \left(
\begin{array}{cc}
-1  & \frac{(r+1)^2}{r}  \\[1mm]
-1  & \frac{(r+1)^2}{r}-1
\end{array}
\right)
\een
as $R\rightarrow +\infty$, and
\be
h_1gh_1^{-1}=\pm \left(
\begin{array}{cc}
 1   & -\frac{(r+1)^2}{r}  \\[1mm]
0 & 1
\end{array}
\right), \hspace*{35mm} \label{h1gh1inv} \\
h_1fg^{-1}h_1^{-1}=\frac{\pm 1}{k(r^2-1)} \left(
\begin{array}{cc}
 -(r+1)\left(k^2-r \right)   & (r+1)^2\left(r  -\frac k r\right)  \\[1mm]
 -r\left(k^2-1\right) &  (r+1)\left(r -k^2\right)
\end{array}
\right) \rightarrow \pm \left(
\begin{array}{cc}
1  & 0  \\
1  & 1
\end{array}
\right) \nn
\ee
as $R\rightarrow +\infty$. Thus, ${l_1}\rightarrow 2\log r$, $l_2\rightarrow 0$. In such a way our once-punctured annulus model $A_1$ becomes the twice-punctured unit disk discussed in Beardon's paper \cite{Beardon2012}. And $t=-\frac{(r+1)^2}{r}$, for $t$ as in \cite[(5)]{Beardon2012}, $r$ as in \eqref{h1gh1inv}. We note that in this case $\cos \theta_1 \rightarrow \tanh \frac{l_1}{4}$ from the first expression of \eqref{main quantity}, then $\tan \theta_1 \rightarrow 1 / \sinh \frac{l_1}{4}$, and for $j=1$, \eqref{comparison ex} becomes
$$\displaystyle \frac{l_1}{\pi}  \leq \lambda_1 \leq \frac{l_1}{\frac{\pi}{2}+\arctan \left( \frac{1}{\sinh ({l_1}/{4})}\right)},$$
which is the same as Theorem 8.2 in \cite{Beardon2012}.

\vspace*{6mm}
{\setlength{\parindent}{0pt}\textbf{Acknowledgement.}} The author would like to thank Prof. Toshiyoki Sugawa for his proposal for this topic, his suggestions and encouragements. The author is also grateful to Prof. Katsuhiko Matsuzaki for his idea in Remark 3, and Prof. Ara Basmajian for his comments on the hyperbolic and complex structures of $A$.
 
\vspace{6mm}
\providecommand{\bysame}{\leavevmode\hbox to3em{\hrulefill}\thinspace}
\providecommand{\MR}{\relax\ifhmode\unskip\space\fi MR }
\providecommand{\MRhref}[2]{%
  \href{http://www.ams.org/mathscinet-getitem?mr=#1}{#2}
}
\providecommand{\href}[2]{#2}

\newpage
\end{document}